\documentclass[11pt,reqno]{amsproc}

\title[]{Pressure, Intermittency, Singularity}

\author{Peter Constantin}
\address{Department of Mathematics, Princeton University, Princeton, NJ 08544}
\email{const@math.princeton.edu}
\usepackage[margin=1in]{geometry}
\usepackage{amssymb,amsmath,amsfonts,amssymb, latexsym, verbatim, esint,amsthm, mathrsfs} 
\usepackage{times}
\usepackage{color}
\usepackage{hyperref}
\newcommand{\pa}{\partial}
\newcommand{\la}{\label}
\newcommand{\fr}{\frac}
\newcommand{\na}{\nabla}
\newcommand{\be}{\begin{equation}}
\newcommand{\ee}{\end{equation}}
\newcommand{\ba}{\begin{array}{l}}
\newcommand{\ea}{\end{array}}
\newcommand{\Rr}{{\mathbb R}}

\newcommand{\beg}{\begin}
\newcommand{\ov}{\overline}

\newcommand{\D}{\Delta}
\renewcommand{\L}{\Lambda}
\renewcommand{\P}{\mathbb P}
\newcommand{\intr}{\int_{\mathbb R^3}}
\date{today}
\begin{document}
\begin{abstract}
We give conditions for regularity of solutions of  three dimensional incompressible Navier-Stokes equations based on the pressure and on structure functions. 
\end{abstract}
\noindent
\maketitle
\begin{center}{\em{On the occasion of the centennial anniversary of O. A. Ladyzhenskaya}}\end{center}
\section{Introduction}
We consider solutions of incompressible Navier-Stokes equations in $\Rr^3$, with smooth and localized initial 
data, and discuss conditions in terms of pressure and structure functions that are easily accessible and guarantee that solutions which are smooth on a time interval $[0, T)$ have smooth (and hence unique) extensions beyond $T$. The literature on regularity issues for Navier-Stokes equations is so extensive that we are not able to give here even the beginning of a survey.  We mention just some minimal references in this short paper, with apologies to the many authors and works we knowingly or unknowingly leave out. \\
We discuss unforced Navier-Stokes equations
\be
\pa_t u + u\cdot\na u -\nu\D u + \na p = 0,
\la{nse}
\ee
with 
\be
\na\cdot u = 0,
\la{divz}
\ee
and 
\be
u(x,0) = u_0.
\la{id}
\ee
The  kinematic viscosity $\nu$ is a strictly positive constant, $u$ is the velocity, $p$ is the pressure. Most of this paper is concerned with solutions in the whole space, but there will be a few instances in which we refer to the bounded domain case. In that case the assumed boundary conditions are homogeneous Dirichlet, 
\be
u_{\left |\right. \,\pa\Omega} = 0.
\la{dirb}   
\ee
 We maintain a sparing notation throughout the paper, omitting arguments and indices as often as we can.\\

 We recall the local existence result for initial data (at time $T_0$) in $V$.  The spaces $H$ (mentioned below) and $V$ are spaces of divergence-free vector fields which are completions of smooth compactly supported divergence free fields in the topologies of $L^2$ and $H^1$. The norm in $V$ is called the enstrophy. In the whole space it corresponds to the $\dot{H}^1$ norm. Initial data with finite enstrophy lead to local strong solutions, that is unique solutions belonging to $L^{\infty}(T_0,  T_0 +\tau ; V)\cap L^2(T_0, T_0+ \tau; H^2\cap V)$ for some $\tau>0$. Strong solutions are $C^{\infty}$ smooth for $t>T_0$ in smooth domains \cite{cfbook}. By "conditions for regularity" for smooth solutions on a time interval $[0,T)$ we mean  conditions which guarantee $u(T)\in V$. These are global regularity conditions.  We note here that we are not talking about $\epsilon$- regularity concepts (\cite{sereginbook}) which are conditions on weak solutions in space-time cylinders, which imply pointwise local regularity inside a smaller cylinder.  When assembled over space time, these conditions lead to partial regularity, and may lead to global regularity if additional assumptions are in place, (for instance a single potential first singularity at one point). In this paper we consider conditions which lead directly to persistence of regularity.\\

There are several well-known conditions for regularity. One of the simplest is
\be
\int_0^T\|u(t)\|_V^4dt  \le M_V <\infty.
\la{vcond}
\ee
From it, we have in a straightforward manner \cite{cfbook} that
\be
\|u(t)\|_V^2 \le \|u(0)\|_V^2\exp{\left(C\nu^{-3} M_V\right)}
\la{quantaMV}
\ee
for all $0\le t\le T$. We adhere to the good practice that arguments of exponentials or logarithms should be nondimensional.  Another easy to prove  explicit condition is based on $\|\na u\|_{L^3}$ (see below, Theorem~\ref{nal3}).

The celebrated  Ladyzhenskaya-Prodi-Serrin conditions \cite{sereginbook}  are
\be
\int_0^T \|u(t)\|_{L^q}^pdt  \le  M_{p,q}<\infty,
\la{pscond}
\ee
with
\be
\fr{2}{p} + \fr{3}{q} = 1,
\la{pqcond}
\ee
and $3<q\le \infty$.   When $q=3$ the condition is
\be
\|u(t)\|_{L^3} \le M_3<\infty, \quad t-a.e. \quad \text{on}\quad [0,T].
\la{l3cond}
\ee
As it is very well-known, the Ladyzhenskaya-Prodi-Serrin conditions imply regularity. The following is the explicit bound on the enstrophy.  
\beg{thm}\la{prodserstr} Let $\Omega$ be a bounded open domain in $\Rr^3$ with smooth boundary, let $q>3$ and let $u$ be a strong solution of the Navier-Stokes equations in $\Omega$ on the interval $[0,T]$. There exists an absolute constant $C$  such that 
\be
\|u(t)\|_V^2 \le \|u(0)\|_V^2\exp{\left[C\nu^{-\fr{q+3}{q-3}}\int_0^t\|u(s)\|_{L^q}^{\fr{2q}{q-3}}ds\right ]}.
\la{quantaV}
\ee
holds for $0\le t\le T$. In particular, if \eqref{pscond} holds then
\be
\|u(t)\|_V^2 \le \|u(0)\|_V^2\exp{\left[C\nu^{-\fr{q+3}{q-3}}M_{p,q}\right ]}.
\la{quantps}
\ee
\end{thm}
\beg{proof}
Here is a brief proof. We recall that the Stokes operator is defined as
\be
Au = -\P \D u,
\la{apu}
\ee
where $\P$ is the Leray  projector on divergence-free vector fields. We recall
\be
\|u\|_{H^2(\Omega)}\le C |Au|_H,
\la{ellipta}
\ee
the fact that
\be
\|u\|_V = \|u\|_{H^1({\Omega})},
\la{uV}
\ee
and the notation
\be
B(u,v) = \P(u\cdot\na v).
\la{buv}
\ee
We start with the enstrophy evolution
\be
\fr{1}{2}\fr{d}{dt}\|u\|^2_V + \nu |Au|^2_H = - (B(u,u), Au)_H \le |B(u,u)|_H|Au|_H
\la{enstri}
\ee
Now, because $\P$ is a projector, it follows that
\be
|B(u,u)|_H \le \|u\cdot\na u\|_{L^2}.\la{bul2}
\ee
A H\"{o}lder inequality with exponents $q, \fr{2q}{q-2}, 2$ yields 
\be
|B(u,u)|_H \le \|u\|_{L^q}\|\na u\|_{L^{\fr{2q}{q-2}}},
\la{buli}
\ee
and because $2<\fr{2q}{q-2}<6$, interpolation yields
\be
\|u\|_{L^q}\|\na u\|_{L^{\fr{2q}{q-2}}} \le \|u\|_{L^q}\|\na u\|_{L^2}^{1-\fr{3}{q}}\|\na u\|_{L^6}^{\fr{3}{q}}
\la{buuinter}
\ee
Using the embedding $H^2(\Omega)\subset W^{1, 6}(\Omega)$, \eqref{ellipta} and \eqref{uV}, we have
\be
|B(u,u)|_H \le  C\|u\|_{L^q}\|u\|_V^{1-\fr{3}{q}}\|Au\|_H^{\fr{3}{q}}
\la{buuinterav}
\ee
Thus, from\eqref{enstri} we have
\be
\fr{1}{2}\fr{d}{dt}\|u\|^2_V + \nu |Au|^2_H \le \|u\|_{L^q} \|u\|_V^{1-\fr{3}{q}}|Au|_H^{1+\fr{3}{q}}
\la{enstro}
\ee
and Young's inequality with exponents $(\fr{1}{2}(1+\fr{3}{q}))^{-1}, (\fr{1}{2}(1-\fr{3}{q})))^{-1}$ yields
\be
\fr{1}{2}\fr{d}{dt}\|u\|^2_V + \nu |Au|^2_H \le \fr{\nu}{2} |Au|^2_H + C\nu^{-\fr{q+3}{q-3}} \|u\|_{L^q}^{\fr{2q}{q-3}} \|u\|_V^{2}.
\la{enstrop}
\ee
The claimed inequality \eqref{quantaV} follows by integrating the ODE inequality \eqref{enstrop}.
\end{proof}
\beg{rem} The same result holds in $\Rr^3$ or $\mathbb T^3$ with the same proof.
\end{rem}
If $q>3$, the bound on the enstrophy is precise and quantitative. In the case $q=3$, in order to have a good quantitative control it is useful to have  a form of finite uniform integrability of $|u(x,t)|^3$.    
This condition is  
\be
\exists \delta>0, \forall t, \forall A, \quad |A|\le \delta \Rightarrow  \int_A|u(x,t)|^3dx \le \left (\fr{\nu}{2C}\right )^3.
\la{uil3}
\ee
In the left hand side, $|A|$ is the Lebesgue measure of $A$. In the right hand side, $\nu$ is the kinematic viscosity and $C$ is  the constant in Morrey's inequality,
\be
\|u\|_{L^6(\Rr^3)} \le C\|u\|_{\dot H^1(\Rr^3)}.
\la{morrey}
\ee
\beg{rem} The condition \eqref{uil3} is uniform in time, but it is much weaker than uniform integrability, because  $\fr{\nu}{2C}$ is fixed.  
\end{rem}

\beg{thm} \la{foiasl3} Let $u$ be a strong solution of the NSE in $\Rr^3$  on $[0,T]$. Assume \eqref{uil3}. Then 
 \be
 \|u(t)\|_{\dot H^1} ^2\le\min\left\{ 
 \ba  \|u_0\|_{\dot H^1}^2\exp{\{t\left(\fr{ \|u_0\|_{L^2}^2}{\delta \nu}\right)\}}, \\
   \|u_0\|_{\dot{H^1}}^2 + \fr{2}{\delta \nu^2}\|u_0\|_{L^2}^4, 
    \ea
   \right.
 \la{ensbdelta}
 \ee
 where $\delta$ is the constant in \eqref{uil3}.
 \end{thm} 
 \beg{rem} The time exponential bound is better than the time independent bound for times shorter than
 $\fr{\delta\nu}{\|u_0\|_{L^2}^2}\log\left(1 + \fr{2\|u_0\|^4_{L^2}}{\delta\nu^2\|u\|^2_{\dot H^1}}\right)$. After that time, the time independent bound is smaller. In either case, the bound \eqref{ensbdelta} implies that the enstrophy is bounded on $[0,T]$, which in turn implies that the solution has a unique strong extension beyond $T$.
 \end{rem}

 \beg{proof}
The proof (based on (\cite{foias})  follows from the enstrophy equation \eqref{enstri} using the fact that 
\be
\left | \{x; |u(x,t)| \ge U\}\right| \le U^{-2}\|u_0\|_{L^2}^2
\la{measbad}
\ee
with the choice of 
\be 
U =\delta^{-\fr{1}{2}}\|u_0\|_{L^2},
\la{Uchoicedelta}
\ee
and estimating the nonlinear term separately in the region where $|u(x,t)| \ge U$ and where $|u(x,t)|\le U$ by
\be
\|u\cdot\na u\|_{L^2} \le \left(\int_{|u||\ge U}|u|^3dx\right)^{\fr{1}{3}}\|\na u\|_{L^6} + U\|\na u\|_{L^2}. 
\la{l3split}
\ee
Then,  using the assumption \eqref{uil3} we obtain
\be
\|u\cdot\na u\|_{L^2}\le \fr{\nu}{2}\|\Delta u\|_{L^2} + U \|\na u\|_{L^2},
\la{unau}
\ee
we absorb the first term in half the dissipation and use a Young inequality in the second term. We end up with ODE inequality
\be
\dot y \le \fr{U^2}{\nu} y
\la{odey}
\ee
for the quantity $y = \|u\|^2_{\dot H^1}$ with $U$ given by \eqref{Uchoicedelta}. The exponential bound in inequality \eqref{ensbdelta} follows from Gronwall, and the time independent bound follows by using $\nu \int_0^t y(s)ds\le 2\|u_0\|^2_{L^2}$.
\end{proof}

\beg{rem} The same result holds in bounded domains $\Omega$ or the periodic case $\mathbb T^3$, with the same proof. We use the enstrophy equation \eqref{enstri}, and the inequality \eqref{bul2}. Then we estimate like in \eqref{l3split} and
note that $\|\na u\|_{L^6} \le C|Au|_H$.
\end{rem}

As the reader may have already noticed, we are interested in explicit conditions, involving constants known a~priori, without the need to sample solutions, and which yield explicit enstrophy bounds.  These conditions are useful if additionally it is true that if they are satisfied uniformly on solutions of approximations that converge only almost everywhere, then the solutions are smooth. We refer to such conditions as "easily accessible". The conditions
\eqref{pscond}, \eqref{uil3} are easily accessible.
\beg{thm} \la{approx} Let $u_n(x,t)$ be an approximation of the NSE solution $u(x,t)$ on the interval $[0,T]$. Assume that
\be
u_n(x,t)\to u(x,t)\quad (x,t)-a.e. \quad \text{on} \quad \Omega\times [0,T].
\la{conv}
\ee

\noindent (i)  If there exists $M_{p,q}<\infty$, such that \eqref{pscond} holds for $u_n$ uniformly for all $n$, then $u$ obeys \eqref{pscond} with the same constant $M_{p,q}$. 

\noindent  (ii) If there exists $M_3$ such that \eqref{l3cond} holds for $u_n$ uniformly for all  $n$, then $u$ obeys \eqref{l3cond} withe the same constant $M_3$.

\noindent (iii) If there exists a constant $\delta>0$ such that \eqref{uil3}  holds for $u_n$ uniformly for all  $n$, then $u$ obeys \eqref{uil3} withe the same constant $\delta$.

\end{thm}

\beg{rem} In the case $(i)$, the solution $u$ obeys the quantitative bound \eqref{quantps}. In the case $(iii)$, the solution $u$ obeys the quantitative bound \eqref{ensbdelta}. In these cases the $H^1$ bound on $u(T)$ is explicit.
\end{rem}
\beg{proof} The proof of $(i), (ii), (iii)$ follows from applications of Fatou's lemma.
\end{proof}

This paper is devoted to conditions based on pressure and on structure functions. There is a good motivation to seek conditions in terms of the pressure. In the absence of the pressure, the Navier-Stokes equations are Burgers equations in 3D and obey a maximum principle. This implies that the velocity is bounded, (if initially so), and the solutions are smooth for all time. The pressure in the Navier-Stokes equations is the only reason the equations are not local and the velocity magnitude is not a~ priori controlled.  Conditions of regularity in terms of the pressure are known in $L^2$ \cite{bergal}  and one sided in $L^{\infty}$ \cite{SerSv}. We present in this paper an $L^{\fr{3}{2}}$ condition, the analogue of the $q=3$ condition \eqref{uil3} expressed in terms of only the pressure (Theorem~\ref{pl32thm}, Theorem ~\ref{presl32quanta}). We also give the analogues of the Ladyzhenskaya-Prodi-Serrin conditions (Theorem ~\ref{prodiserrinpl3}).\\

The motivation to express conditions for regularity in terms of structure functions comes from experimental, numerical and theoretical turbulence studies. Structure functions are averages of moments of velocity increments. They obey remarkable and robust statistical relations. The relations need interpretation and then the may serve as reasonable hypothesis for the solutions of Navier-Stokes equations. One of the more widely verified relations is the "four-fifths" law \cite{frisch}
\be
\langle (\delta^{\parallel}_{\ell} u)^3\rangle = -\fr{4}{5}\epsilon |\ell |
\la{45}
\ee
where $\delta^{\parallel}_{\ell} (u) = (u(x+\ell) - u(x))\cdot \fr{\ell}{|\ell|}$ is the longitudinal velocity increment and  $\epsilon = -\langle \fr{dE}{dt}\rangle$ is the rate of dissipation of energy, which in the case of unforced NSE
equals $\nu \langle |\na u |^2\rangle$.  The four-fifths law is shown to hold for homogeneous and isotropic turbulence in the limit of time to infinity, followed by Reynolds number to infinity, followed by $\ell \to 0$, in this order. The Navier-Stokes solutions are assumed to be smooth. The braces $\langle\cdot\rangle$ are expectations (ensemble average). Long time and space averages are usually assumed to realize them, and in numerical experiments these averages are used. The assumption of finite positive $\epsilon$ is also made, in the limit of time to infinity, followed by Reynolds number to infinity, in this order. The Reynolds number is defined as 
\be
Re=\fr{UL}{\nu}  
\la{rey}
\ee
where $U$ is a velocity scale, and $L$ is a length scale.  The classical K'41 Kolmogorov theory proposes scaling exponents $\zeta_p = \fr{p}{3}$ for structure functions
\be
\langle (\delta^{\parallel}_{\ell} u)^p\rangle = C_p (\epsilon |\ell|)^{\zeta_p}.
\la{spzetap}
\ee
These relations are expected to hold in a range of scales, $ |\ell |\in (\eta, L)$ where $\eta$ is the Kolmogorov dissipation scale,
\be
\eta = \left(\fr{\nu^3}{\epsilon}\right)^{\fr{1}{4}}
\la{eta}
\ee
which is determined by the kinematic viscosity and energy dissipation rate, alone. Below the Kolmogorov dissipation scale, it is assumed that viscous effects dominate, with smooth behavior. The length scale $L$ is the integral scale of turbulence. 
Turbulence findings are average statements, they refer to typical long time behavior, and are asymptotic in Reynolds number. Interpreting them for the initial value problem for Navier-Stokes equations is challenging. It is however reasonable to expect that there are many solutions which give statistical weight to the turbulence laws and have properties that are consistent with them.

We give quantitative conditions in Theorem~\ref{S2reg}, Theorem~\ref{S2quanta}. They involve a cutoff scale $r=r(t)$. We modify the structure function $S_2(x, r)$ (see \eqref{S2} below) to take into account a possibly non-universal  viscous regularization below $r$. The regularity condition \eqref{S2cond} requires $\int_0^T r(t)^{-4}dt <\infty$, a condition satisfied by the Kolmogorov length $r= \eta$. The condition requires in addition the smallness of $\int_A S_2^{\fr{3}{2}}(x,r)dx \le \left(\fr{\nu}{C}\right)^3$ on sets of small enough measure, $|A|\le \delta$.

Modern theories modify the scaling $\zeta_p = \fr{p}{3}$  in \eqref{spzetap} of the K'41 theory, reflecting  experimental and numerical observation of intermittency. The turbulent signal is intermittent, that is, regions of high gradients of velocity are found to be sparse in both time and space. The connection between intermittency and regularity was explored in several mathematical works, (see for instance \cite{gruj} and references therein) where assumptions of  sparse behavior in physical space are used to deduce improved conditional regularity. In a different setting \cite{chesk}, multifractal scaling exponents were connected to ratios of volume averages, and conditions for regularity were given on the basis of intermittency dimension.

In terms of the exponents, it is found numerically (see for instance \cite{sreeni}) that $\zeta_2>\fr{2}{3}$ and, while $\zeta_3$ remains close to 1, $\zeta_p$ become smaller than $p/3$ for large $p$, and perhaps even tends asymptotically to a constant, suggesting depletion of regularity, significantly below H\"{o}lder. We express conditions of regularity in terms of a Dini modulus of continuity which give regularity if logarithmic scaling is assumed (Theorem~\ref{dini}).  In Section~\ref{multi} we consider a multifractal scenario where regularity still persists. In Section \ref{time} we give a condition for regularity which requires small increments of velocity only in time dependent regions of high velocity and high gradients.

The proofs are based on observations concerning the pressure. 

\section{The pressure}
We consider solutions of
\be
-\Delta p = \na\cdot(u\cdot\na u)
\la{peq}
\ee
in $\Omega\subset \Rr^3$, where $u$ is divergence-free and sufficiently regular.  We recall representation results from \cite{pres}. We use the notation
\be
{\ov{f}}(x,r) = \fr{1}{4\pi r^2}\int_{|x-y|=r}f(y)dS(y) = \fint_{|\xi| =1} f(x+r\xi)dS(\xi)
\la{ovf}
\ee
where $dS$ is surface area and $\fint$ denotes the integral normalized by the area of the region of  integration, which in the above case is $4\pi$.
We denote 
\be
\sigma_{ij}(\widehat{y-x}) = \frac{3(y_i-x_i)(y_j-x_j)}{|y-x|^2}-\delta_{ij}
\la{sigma}
\ee
where
\be
{\widehat{{y-x}}} = \fr{y-x}{|y-x|}.
\la{hatyx}
\ee
The following lemma was proved  in \cite{pres}.

\beg{lemma} \la{lem2} Let $x\in\Omega\subset \Rr^3$, let $0<r<{\mbox{dist}}(x,\partial \Omega)$, and let $p$ solve (\ref{peq}) with divergence-free $u\in C^{2}(\Omega)^3$. Let $v\in \Rr^3$. Then
\be
\ba
p(x)  =   {\ov{p}}(x,r)  - \fr{1}{3} |u(x)-v|^2 + \fint_{|\xi|=1}\left |\xi\cdot(u(x+ r\xi)- v)\right |^2dS(\xi)  \\ +
P.V. \int _0^r\fr{d\rho}{\rho}\fint_{|\xi|=1}\sigma_{ij}(\xi)(u_i(x+\rho\xi)-v_i)(u_j(x + \rho\xi)- v_j)dS(\xi).
\ea
\la{pv}
\ee 
\end{lemma}
All terms in the right hand side of \eqref{pv} are determined solely by information in the ball of radius $r$ about $x$. We denote the singular integral 
\be
K(x,r) = P.V. \fr{1}{4\pi}\int_{|y-x|\le r} \fr{\sigma_{ij}\left(\widehat{y-x}\right)}{|y-x|^3}(u_i(y)-v_i)(u_j(y)- v_j)dy.
\la{K}
\ee
Thus, \eqref{pv} reads
\be
p(x)  =  {\ov{p}}(x,r) - \fr{1}{3} |u(x)-v)|^2  +  \fint_{|\xi|=1}\left |\xi\cdot(u(x+ r\xi)- v)\right |^2dS(\xi) + K(x,r).
\la{pbs2}
\ee
Evidently, $K$ depends on the choice of the vector $v$. In applications we want to be able to choose $v$ appropriately.  For divergence free functions $u$, it holds that 
\be
\fint_{|\xi |= 1}\xi_i\left(\xi\cdot u(x+ r\xi)\right )dS(\xi) + \fr{1}{4\pi}PV\int_{B(x,r)}\fr{\sigma_{ij}(\widehat{x-y})}{|x-y|^3}u_j(y)dy = \fr{1}{3}u_i(x)
\la{usig}.
\ee
This follows from 
\[
\fr{1}{4\pi}\int_{B(x,r)}\fr{y_i-x_i}{|y-x|^3} (\na\cdot u)(y)dy = 0
\]
by integration by parts.  We also note that 
\be
\fint_{|\xi |=1}(\xi\cdot v)^2dS(\xi) = \fr{1}{3}|v|^2.
\la{avsq}
\ee
This is true for any $v$ that does not depends on $\xi$. Therefore, the representation \eqref{pv} is valid even if $v$ is a function of $x$ and $r$ (but not $\xi$). Indeed, this follows by opening brackets in the right hand side of \eqref{pv}, using \eqref{usig} and
\eqref{avsq}, and identifying what remains as \eqref{pv} for $v=0$, which was proved independently in \cite{pres}.

We average \eqref{pbs2} $\fr{1}{R}\int_R^{2R}dr$  and obtain the representation \cite{pres}
\beg{thm}\la{rep} Let $p$ solve \eqref{peq} with divergence-free $u\in C^2(\Omega)^3$. Let $x\in \Omega\subset \Rr^3$,  $v\in \Rr^3$ and let $0<r<\fr{1}{2}{dist}(x,\pa\Omega)$. Then, 
\be
p(x) = \beta(x,r) + \pi(x,r) - \fr{1}{3r}\int_r^{2r}|u(x)-v|^2d\rho
\la{pbetapi}
\ee
with
\be
\beta(x,r) = \fr{1}{r}\int_r^{2r}\ov{p}(x,\rho)d\rho
\la{beta}
\ee
and
\be
\pi(x,r) = \fr{1}{r}\int_r^{2r}\left[K(x,\rho) + \fint_{|\xi| =1} |\xi\cdot (u(x+\rho\xi)-v)|^2dS(\xi)\right]d\rho
\la{pixrk}
\ee
The explicit expression for $\pi(x,r)$  is
\be
\ba 
\pi(x,r) =
 P.V.\fr{1}{4\pi}\int_{|x-y|\le 2r} w\left(\fr{|y-x|}{r}\right)\fr{\sigma_{ij}(\widehat{y-x})}{|y-x|^3}(u_i(y)-v_i)(u_j(y)- v_j)dy \\
 + \fr{1}{4\pi r}\int_{r\le |y-x|\le 2r}\fr{1}{|y-x|^2}\left (\fr{y-x}{|y-x|}\cdot (u(y)-v)\right)^2dy,
\ea
\la{pixr}
\ee
where the weight $w$ is given by
\be
w(\lambda) =
\left\{
\ba
1,\quad \quad{\rm{if}}\; \; 0\le \lambda \le 1,\\
2 -\lambda \quad{\rm{if}}\;\; 1\le \lambda\le 2,\\
0 \quad\quad {\rm if}\;\; \lambda \ge 2
\ea
\right.
\la{weight}
\ee  
\end{thm}
We recall  bounds on $\beta$  and $\pi$ (Propositions 2 and 3, \cite{pres})
\beg{prop}\la{betabound}
There exists an absolute constant $C$ such that, for any $r>0$, 
\be
\|\na\beta (\cdot, r)\|_{L^2(\Rr^3)} \le C r^{-1}\|u\|_{L^4(\Rr^3)}^2
\la{nabetab}
\ee
and 
\be
\|\beta (\cdot, r)\|_{L^{\infty}}\le Cr^{-2}\|\na u\|_{L^2}\|u\|_{L^2}
\la{betifty}
\ee
hold. 
Moreover, for any $1<q<\infty$, there exists $C_q$ independent of $r$, such that
\be
\|\beta(\cdot, r)\|_{L^q(\Rr^3)} \le C_q \|u\|_{L^{2q}(\Rr^3)}^2
\la{lqbounds}
\ee
holds.
\end{prop}
We choose now  $v = u(x)$. Then from \eqref{lqbounds} and the corresponding bound for $p$ it follows also that
\be
\|\pi(\cdot, r)\|_{L^q(\Rr^3)} \le C_q \|u\|_{L^{2q}(\Rr^3)}^2.
\la{pilqbounds}
\ee
We denote by
\be
S_2(x,r) = \fr{1}{4\pi}\int_{|y|\le 2r} \fr{1}{|y|^3}|u(x+y)-u(x)|^2dy
\la{S2}
\ee.
We note that
\be
|\pi(x,r)|\le 2 S_2(x, r)
\la{pis2b}
\ee
follows from \eqref{pixr}.  For any measurable set $A\subset\Omega$ with $dist(A, \pa\Omega)>2r$ we have
\be
\int_A S_2(x,r)^q \le C_q r^{2q}\int_{A+ rB(0,1)} |\na u(x)|^{2q}dx
\la{suna}
\ee
for any $q\ge 1$ where $B(0,1)$  is the unit ball in $\Rr^3$. This follows in straightforward manner by writing the integral in \eqref{S2} in polar coordinates with $y=\rho\xi$, $\rho = |y|$, $\xi =\widehat{y}$, expressing
$u(x+y)- u(x) = \int_0^1 \fr{d}{d\lambda} u(x +\lambda\rho\xi)d\lambda$, and  using Schwarz and H\"{o}lder inequalities.

\beg{rem} The bounds \eqref{lqbounds} and \eqref{pilqbounds} for $\beta$ and $\pi$ are valid in bounded domains with smooth boundary if we add bounds
for $\|p\|_{L^q(\pa\Omega)}$, see Lemma 2 in \cite{esc}. Once these are obtained the rest of the bounds which are local bounds are valid.
\end{rem}

\section{Conditional regularity}
\beg{thm} \la{nal3} Assume 
\be
\int_0^T \|\na u(t)\|_{L^3}^2 dt  = N <\infty.
\la{A}
\ee
There exists an absolute constant $C$ such that
\be
\sup_{0\le t \le T} \|\na u(t)\|_{L^2}^2 \le  \|\na u_0\|_{L^2}^2 \exp{\fr{CN}{\nu}}
\la{enstrophyboundA}
\ee
holds.
\end{thm}
\beg{proof}
\[
 \fr{1}{2}\fr{d}{dt} \intr |\na u |^2 dx + \nu\int |\D u|^2 dx = \intr (u\cdot\na u)\cdot\D u dx\le \|u\|_{L^6}\|\na u\|_{L^3}\|\D u\|_{L^2}
\]
followed by Schwartz and Morrey inerqualities results in
\[
 \fr{1}{2}\fr{d}{dt} \|\na u (t)\|_{L^2}^2 \le \fr{C}{\nu}  \|\na u(t)\|_{L^3}^2 \|\na u(t)\|_{L^2}^2,
 \]
 and the proof is finished by Gronwall.
 \end{proof}
 Let 
 \be
 \L = \left(-\D\right)^{\fr{1}{2}}.
\la{Lambda}
\ee
 \beg{thm} \la{na32l2}
 Assume 
\be
\int_0^T \|\L^{\fr{3}{2}}u(t)\|_{L^2}^2 dt  = M<\infty.
\la{M}
\ee
There exists an absolute constant $C$ such that
\be
\sup_{0\le t \le T} \|\na u(t)\|_{L^2}^2 \le  \|\na u_0\|_{L^2}^2 \exp{\fr{CM}{\nu}}
\la{enstrophyboundM}
\ee
holds.
\end{thm}
\beg{proof} This is a consequence of \eqref{enstrophyboundA} and of the inequality $N\le CM$ which follows from the  familiar \cite{stein} Riesz potential inequality
\be
\|f\|_{L^3} \le C \|\L^{\fr{1}{2}}f\|_{L^2}
\la{homoh1/2}
\ee
with $f = \pa_i u_j$, for $i,j =1, \dots 3$.
\end{proof} 

 \beg{rem} The finiteness \eqref{enstrophyboundA} is a stronger regularity condition (i.e., more general) than the finiteness \eqref{enstrophyboundM}. Both conditions are analogues of the Ladyzhenskaya-Prodi-Serrin  condition \eqref{pscond} for $q=\infty $.  The finiteness  \eqref{enstrophyboundA} and \eqref{pscond} for $q=\infty $ are logically independent of each other.
 \end{rem}
 
 \beg{thm} 
There exists an absolute constant $C$ such that
\be
\sup_{0\le t \le T} \|\L^{\fr{1}{2}} u(t)\|_{L^2} \le  \|\L^{\fr{1}{2}} u_0\|_{L^2} +  C\int_0^T  \|\L^{\fr{3}{2}}u(t)\|_{L^2}^2 dt  
\la{homoh12bound}
\ee
holds for smooth solutions of Navier-Stokes equations.
\end{thm}
 \beg{proof}
 We take the scalar product of the Navier-Stokes equation with $\L u(t)$ and integrate. We obtain
 \be
 \fr{1}{2}\fr{d}{dt} \|\L^{\fr{1}{2}} u\|_{L^2}^2 + \nu \|\L^{\fr{3}{2}} u\|_{L^2}^2 \le \left | \int (u\cdot\na u)\cdot \L u dx\right |
 \le C\|u\|_{L^3}\|\na u\|_{L^3}\|\L u\|_{L^3} \le C \|\L^{\fr{1}{2}} u\|_{L^2}\|\L^{\fr{3}{2}} u\|_{L^3}^2,
 \ee
 where we used \eqref{homoh1/2}. We divide by $\|\L^{\fr{1}{2}} u\|_{L^2}$ and integrate in time.
 \end{proof}
 \beg{rem}
 In view of \eqref{homoh1/2}, the inequality \eqref{homoh12bound} implies the boundedness of the $L^3$ norm. Hence, via the Navier-Stokes equation, the finiteness \eqref{M} implies  \eqref{l3cond}.
 The main virtue of \eqref{homoh12bound} is that it is viscosity-independent, and is valid for Euler equations as well.
 \end{rem}

 The analogue of the $L^3$-based condition \eqref{l3cond} in terms of only the pressure is the following.
\beg{thm}\la{pl32thm}
There exists an absolute constant $C$, such that, if $p= R_iR_j(u_iu_j)$ satisfies the finite uniform integrability condition
\be
\exists \delta>0,\forall t, \forall A \quad  |A|\le \delta \Rightarrow \int_A|p(x,t)|^{\fr{3}{2}}dx \le \left(\fr{\nu}{C}\right)^3
\la{unifint}
\ee
on $ [0,T]$, then $u\in L^{\infty}(0,T; L^3(\Rr^3))$ with explicit bounds depending only on $ \delta, \nu, T, \|u_0\|_{L^3}, \|u_0\|_{L^2}$.
\end{thm}
\beg{proof} We take the evolution of the $L^3$  norm of $u$:
\be
\fr{1}{3} \fr{d}{dt}\intr |u|^3 dx + \nu\intr |u| [|\na u|^2 + |\na |u||^2]dx = -\intr |u|(u\cdot\na p) dx = \intr p(u\cdot\na |u|)dx
\la{l3ev}
\ee
Let $U$ be a large positive number and $\phi$ a positive smooth function of one variable which is compactly supported
in $[0,2]$, satisfies  $0\le \phi\le 1$ and identically equals $1$ on $[0,1]$. We split the RHS of \eqref{l3ev} ,
\[
\intr p(u\cdot\na |u|) dx = \intr \phi\left(\fr{|u|}{U}\right)p(u\cdot\na |u|)dx + \intr\left (1-\phi\left(\fr{|u|}{U}\right)\right)p(u\cdot\na |u|)dx. 
\]
We estimate the first term, using that on the support of $\phi$ we have $|u|\le 2U$,
\be
\left |\intr \phi\left(\fr{|u|}{U}\right)p(u\cdot\na |u|)dx\right| \le (2U)^{\fr{1}{2}}\sqrt{D}\sqrt{\intr p^2dx} 
\la{oneone}
\ee
where
\be
D = \intr |u| |\na u|^2dx
\la{D}
\ee
and then using the boundedness of Riesz transforms in $L^p$ spaces and then interpolating $L^4$ between $L^3$ and $L^9$,  we have
\be
\intr p^2\le C \intr |u|^4dx \le C \|u\|_{L^3}^{\fr{5}{2}}\|u\|_{L^9}^{\fr{3}{2}}.
\la{pl2interp}
\ee
Now we use the fact that there exists a constant $C$ such that
\be
D \ge C\|u\|_{L^9}^3.
\la{dl9}
\ee
This fact follows from Morrey's inequality $\|\na f\|_{L^2}\ge C\|f\|_{L^6}$ applied for
with  $f= |u|^{\fr{3}{2}}$.
Thus, from \eqref{oneone} and \eqref{dl9} we have
\be
\left |\intr \phi\left(\fr{|u|}{U}\right)p(u\cdot\na |u|)dx\right | \le CU^{\fr{1}{2}}D^{\fr{3}{4}}\|u\|_{L^3}^{\fr{5}{4}}.
\la{onetwo}
\ee
With Young's inequality we have
\be
\left |\intr \phi\left(\fr{|u|}{U}\right)p(u\cdot\na |u|)dx\right | \le \fr{\nu}{2} D + C U^2\nu^{-3} \|u\|_{L^3}^5.
\la{oneall}
\ee
We know from the energy inequality, interpolation $L^2-L^6$ and Morrey's inequality that  
\be
\int_0^T \|u\|_{L^3}^4 dt \le C\|u_0\|_{L^2}^2\int_0^T \|\na u\|_{L^2}^2 dt \le C\nu^{-1}\|u_0\|^4_{L^2},
\la{intl34}
\ee
and thus the factor $CU^2 \nu^{-3}\|u\|_{L^3}^2$ multiplying the $\|u\|_{L^3}^3$ is time integrable.
The second term is estimated as
\be
\left | \intr\left (1-\phi\left(\fr{|u|}{U}\right)\right)p(u\cdot\na |u|)dx\right|\le C\sqrt{D}\|u\|_{L^9}^{\fr{1}{2}}\|((1-\phi)p\|_{L^{\fr{9}{4}}}
 \la{twoone}
 \ee 
 by using H\"{o}lder with exponents $\fr{9}{4},18, 2 $. One more interpolation, 
 \[
 \|(1-\phi)p\|_{L^{\fr{9}{4}}}\le \|(1-\phi)p\|_{L^{\fr{9}{2}}}^{\fr{1}{2}}\|(1-\phi) p\|_{L^{\fr{3}{2}}}^{\fr{1}{2}}
 \]
 and the inequality based on the fact that Riesz transforms are bounded in $L^p$
 \[
 \|(1-\phi)p\|_{L^\fr{9}{2}}\le C \|u\|_{L^9}^2,
 \]
 together with \eqref{dl9}
 yield from \eqref{twoone}
 \be
 \intr\left (1-\phi\left(\fr{|u|}{U}\right)\right)p(u\cdot\na |u|)dx\le C{D}\|((1-\phi)p\|_{L^{\fr{3}{2}}}^{\fr{1}{2}}.
 \la{twoall}
\ee
Now the support of $1-\phi\left(\fr{|u(x,t)|}{U}\right)$ is included in the set
\be
B_U(t) = \{x\left |\right.\; |u(x,t)| \ge U\}
\la{bu}
\ee
which has uniformly small Lebesgue measure
\be
|B_U(t)| \le U^{-2}\|u_0\|^2_{L^2}
\la{bumeas}
\ee
and, by assumption, the function $x\mapsto |p(x,t)|^{\fr{3}{2}}$ satisfies \eqref{unifint}, so that
\be
\|((1-\phi)p\|_{L^{\fr{3}{2}}}^{\fr{1}{2}}\le \fr{\nu}{2C}
\la{smallp}
\ee
holds uniformly for $t\in [0,T]$, if $U$ is chosen large enough.
\end{proof}
\beg{rem}
The condition \eqref{unifint} is weaker than uniform integrability, because $\fr{\nu}{C}$ is fixed. The condition holds if $|p(x,t)|^{\fr{3}{2}}$ is uniformly integrable  in $x$ on $[0,T]$. In particular, it holds if
\be
p\in C(0,T; L^{\fr{3}{2}}(\Rr^3)),
\la{pl32}
\ee
or if $p$ is piece-wise continuous on $[0,T]$ with values in $L^{\fr{3}{2}}$, because in these cases the curve $t\mapsto p$ belongs to a compact subset of $L^{\fr{3}{2}}$ and therefore  $|p(x,t)|^{\fr{3}{2}}$ is uniformly integrable in $x$ on $[0,T]$. The condition also holds if $|p(x, t)| \le f(x)$, $x$-a.e.,  with $f\in L^{\fr{3}{2}}(\Rr^3)$ time independent, because this again implies uniform integrabilty. 
\end{rem}
\beg{rem}
As it is seen in the proof above, the condition \eqref{unifint}  is not applied on just any set, but rather on a set of interest, namely the set where the absolute magnitude of velocity exceeds a fixed large threshhold,  \eqref{bu}. It is  easy to see that this set can be replaced by a smaller, and even more interesting set
\be
B_{U,G}(t)  = \{x\left | \right.\;  |u(x,t)| \ge U, \; \text{and}\; |\na u(x,t)| \ge G\}.
\la{bug}
\ee
The proof of this fact is similar to the proof above, using two cutoffs, and showing that the regions $|u|\le U$ and, separately $|\na u| \le G$ lead each to a priori bounds on the size of the $L^3$ norm of $u$, leaving only the contributions from $B_{U,G}(t)$ to require control.  
\end{rem}
 The following result shows that the condition \eqref{unifint} leads to easily accessible bounds.
\beg{thm}\la{presl32quanta} Let $r\ge 4$. There exists a constant $C=C_r$ such that if \eqref{unifint} holds on $[0,T]$
then
\be
u\in L^{\infty}(0,T; L^r(\Rr^3))
\la{urbl32pr}
\ee
holds. More precisely, we have the single exponential bound
\be
\|u(\cdot,t)\|_{L^r} \le \|u_0\|_{L^r} \exp{\left (\fr{C t\|u_0\|^2_{L^2} }{\nu\delta}\right)}
\la{lrbound}
\ee
where $\delta$ is the constant in \eqref{unifint}.
\end{thm}
\beg{proof}  The evolution of the $L^r$ norm is given by
\be
\fr{1}{r}\fr{d}{dt}\intr |u|^r dx + \nu\intr [|\na u|^2 |u|^{r-2}  + (r-2) |\na |u| |^2 |u|^{r-2}]dx  = \intr pu\cdot\na |u|^{r-2}dx.
\la{lrev}
\ee  
We let $U$ be a large positive number and $\phi$ a positive smooth function of one variable which is compactly supported
in $[0,2]$, satisfies  $0\le \phi\le 1$ and identically equals $1$ on $[0,1]$. We split the RHS of \eqref{l3ev} . We split the right hand side of \eqref{lrev},
\[
\intr p(u\cdot\na |u|^{r-2}) dx = \intr \phi\left(\fr{|u|}{U}\right)p(u\cdot\na |u|^{r-2})dx + \intr\left (1-\phi\left(\fr{|u|}{U}\right)\right)p(u\cdot\na |u|^{r-2})dx. 
\]
We estimate the first term using H\"{o}lder inequalities with exponents $\fr{r}{2}$ for $p$, $2$ for a term involving the gradient, $|u|^{\fr{r-2}{2}}|\na u|$,  and $\fr{2r}{r-4}$ for the term $|u|^{\fr{r-2}{2}}$, taking advantage of the fact that on the support of $\phi$ we have $|u|\le 2U$, and using the boundedness of Riesz transforms in $L^p$ spaces. We deduce that the first term is bounded  by
\be
\left |\intr \phi\left(\fr{|u|}{U}\right)p(u\cdot\na |u|^{r-2})dx \right| \le C U\|u\|_{L^r}^{\fr{r}{2}}D^{\fr{1}{2}},
\la{rle4}
\ee
where
\be
D = \intr |u|^{r-2} |\na u|^2dx.
\la{Dlr}
\ee
We used, in view of $\fr{r(r-2)}{r-4} = r + \fr{2r}{r-4}$, that
\be
\left(\int_{|u|\le U}|u|^{\fr{r(r-2)}{r-4}}dx\right)^{\fr{r-4}{2r}} \le U\|u\|_{L^r}^{\fr{r}{2}-2}.
\la{uUq}
\ee
 The bound \eqref{rle4} is valid for $r=4$ as well, we just take $|u|\le U$ outside the integral and use $L^2-L^2$ bounds.
Hiding $\sqrt{D}$  in $\fr{1}{2}\nu D$, we see that the inequality \eqref{rle4}  leads to an exponential growth 
 \be
 \|u(\cdot,t)\|_{L^r} \le \|u_0\|_{L^r} e^{C\nu^{-1}\int_0^t U^2ds}
\la{expolr}
\ee
if the second term does not contribute to growth,   The second term is bounded using
 \be
D \ge C\|u\|_{L^{3r}}^r.
\la{dl3r}
\ee
We first bound
\be
\ba
\left |  \intr\left (1-\phi\left(\fr{|u|}{U}\right)\right)p(u\cdot\na |u|^{r-2})dx\right | \le \|\left (1-\phi\left(\fr{|u|}{U}\right)\right)p\|_{L^\fr{3r}{r+1}}\|u\|_{L^{3r}}^{\fr{r-2}{2}} D^{\fr{1}{2}}\\
 \le C\|\left (1-\phi\left(\fr{|u|}{U}\right)\right)p\|_{L^\fr{3r}{r+1}} D^{1-\fr{1}{r}},
\ea
\la{twooner}
\ee
and then use 
\be
\ba
\|\left (1-\phi\left(\fr{|u|}{U}\right)\right)p\|_{L^\fr{3r}{r+1}}\le \|\left (1-\phi\left(\fr{|u|}{U}\right)\right)p\|_{L^\fr{3}{2}}^{\fr{1}{2}}
\|\left (1-\phi\left(\fr{|u|}{U}\right)\right)p\|_{L^\fr{3r}{2}}^{\fr{1}{2}}\\
\le C\|\left (1-\phi\left(\fr{|u|}{U}\right)\right)p\|_{L^\fr{3}{2}}^{\fr{1}{2}}D^{\fr{1}{r}}
\ea
\ee
to deduce
\be
\left |  \intr\left (1-\phi\left(\fr{|u|}{U}\right)\right)p(u\cdot\na |u|^{r-2})dx\right |\le C\|\left (1-\phi\left(\fr{|u|}{U}\right)\right)p\|_{L^\fr{3}{2}}^{\fr{1}{2}} D.
\la{twotwor}
\ee
The proof is completed by the assumption of finite uniform integrability \eqref{unifint} which shows this term to be absorbed in the remaining dissipative term $\fr{1}{2}\nu D$.
\end{proof}
\beg{rem}\la{quantaprem} The result above can be used  together with  Theorem \eqref{prodserstr} to give a quantitative bound depending on $\delta$ on the supremum in time of the enstrophy.
\end{rem}
The following result is the analogue of the Ladyzhenskaya-Prodi-Serrin condition in terms of the pressure.
\beg{thm}\la{prodiserrinpl3} Let $p=R_iR_j(u_iu_j)$. Assume that there exists $q>\fr{3}{2}$ such that
 \be
 \int_0^T \|p(t)\|_{L^q(\Rr^3)}^{\fr{2q}{2q-3}}dt <\infty
 \la{plq}
 \ee
Then $u\in L^{\infty}(0,T; L^3(\Rr^3))$ obey $u\in L^{\infty}(0,T; L^3(\Rr^3))$ with explicit bounds depending in addition to \eqref{plq} only on $\nu, T, \|u_0\|_{L^3}, \|u_0\|_{L^2}$.
\end{thm}
\beg{proof} The proof follows from the evolution of the $L^3$ norm \eqref{l3ev} by estimating as in \eqref{twoone}
using H\"{o}lder with exponents $9/4, 18, 2$, 
\be
\left | \intr p(u\cdot\na |u|)dx\right |\le C\sqrt{D}\|u\|_{L^9}^{\fr{1}{2}}\|p\|_{L^{\fr{9}{4}}} \le CD^{\fr{2}{3}}\|p\|_{L^{\fr{9}{4}}}
 \la{puna}
 \ee
 where we used also \eqref{dl9}.  We distinguish three ranges of $q$.
 
 When $q\le \fr{9}{4}$ we interpolate
 \be
 \|p\|_{L^{\fr{9}{4}}} \le \|p\|_{L^q}^{\fr{2q}{9-2q}}\|p\|_{L^{\fr{9}{2}}}^{\fr{9-4q}{9-2q}}.
 \la{interpoq}
 \ee
We use $\|p\|_{L^{\fr{9}{2}}}\le C D^{\fr{2}{3}}$ which folllows from \eqref{dl9} and the boundedness of Riesz transforms in $L^p$ spaces, to deduce
\be
\left | \intr p(u\cdot\na |u|)dx\right |\le C D^{1-\alpha} \|p\|_{L^q}^{\fr{2q}{9-2q}}
\la{pshtaim}
\ee
with $\alpha = \fr{2q-3}{(9-2q)}$. From Young's inequality, the boundedness of the $L^3$ norm of $u$ follows if we know that $\int_0^T \|p\|_{L^q}^{\fr{2q}{2q-3}}dt $ is finite, which was assumed in \eqref{plq}.

When $q\in [\fr{9}{4}, \fr{9}{2}]$, we use Young's inequality in \eqref{puna} and deduce that we need to estimate
$\int_0^T\|p\|_{L^{\fr{9}{4}}}^3dt$. We interpolate
\be
\|p\|_{L^{\fr{9}{4}}}\le \|q\|_{L^q}^{\fr{2q}{3(2q-3)}}\|p\|_{L^{\fr{3}{2}}}^{\fr{4q-9}{3(2q-3)}}\le  
C\|q\|_{L^q}^{\fr{2q}{3(2q-3)}}\|u\|_{L^{3}}^{\fr{8q-18}{3(2q-3)}}
\la{qinterp}
\ee
and thus
\be
\|q\|_{L^{\fr{9}{4}}}^3 \le \|p\|_{L^q}^{\fr{2q}{2q-3}}\|u\|_{L^3}^{\alpha}
\la{q94}
\ee
with $\alpha = \fr{8q-18}{2q-3}$.  If $q\le \fr{9}{2}$ we have $\alpha \le 3$, and the condition \eqref{plq} ensures that $\|u\|_{L^3}$ remains bounded.

When $q\ge \fr{9}{2}$ we estimate
\be
\left | \intr p(u\cdot\na |u|)dx\right |\le C\sqrt{D}\|u\|_{L^3}^{\fr{1}{2}}\|p\|_{L^{3}} 
 \la{pul3}
 \ee
using H\"{o}lder with exponents $2,3,6$. Using Young's inequality we need to consider the effect of the quantity
$\|u\|_{L^3}\|p\|_{L^3}^2$.  We interpolate
\be
\|p\|_{L^3} \le \|p\|_{L^q}^{\fr{q}{2q-3}}\|p\|_{L^{\fr{3}{2}}}^{\fr{q-3}{2q-3}} \le \|p\|_{L^q}^{\fr{q}{2q-3}}\|u\|_{L^{3}}^{\fr{2(q-3)}{2q-3}} \la{pint},
\ee
and it follows that 
\be
\|u\|_{L^3}\|p\|_{L^3}^2 \le C \|p\|_{L^q}^{\fr{2q}{2q-3}}\|u\|_{L^3}^{\fr{6q-15}{2q-3}}
\la{p3f}
\ee
Because $\fr{6q-15}{2q-3} <3$, the condition \eqref{plq} implies a uniform bound on $\|u\|_{L^3}$ on $[0,T]$ in this last case.
\end{proof}

\section{Structure function}
We assume 
\be
\ba 
\exists r(t), \int_0^T r(t)^{-4}dt <\infty, \exists \delta >0, \forall A, \forall  t\in [0,T] \\
|A| \le \delta \Rightarrow \int_A S_2(x,2r(t))^{\fr{3}{2}}dx \le \left(\fr{\nu}{C}\right)^3
\ea
\la{S2cond}
\ee
holds where $S_2(x,r)$  is given in \eqref{S2}.
\beg{thm}\la{S2reg} There exists an absolute constant $C$, such that, for any $T>0$, if a strong solution $u$ 
of NSE satisfies \eqref{S2cond} for all $0\le t<T$, then $u\in L^{\infty}(0, T; L^3(\Rr^3))$  with explicit bounds depending only on $\nu, T, \|u_0\|_{L^3}, \|u_0\|_{L^2}$ and the assumed $\delta>0$, $\int_0^T r(t)^{-4}dt$.
\end{thm}
\beg{proof}
We use \eqref{l3ev} and decompose the pressure $p = \pi + \beta$ as in  \eqref{pbetapi}, at each time $t$, with the choice $r=r(t)$. We bound the term
\be
\left |\intr |u|(u\cdot \na\beta) dx\right| \le Cr^{-1}\|u\|_{L^4}^4
\la{betaint}
\ee
using \eqref{nabetab}, and then by interpolation (used also in \eqref{pl2interp}) we obtain
\be 
\left |\intr |u|(u\cdot \na\beta) dx\right| \le Cr^{-1}\|u\|_{L^3}^{\fr{5}{2}}\|u\|_{L^9}^{\fr{3}{2}}.
\la{betaintb}
\ee
The dissipation $D$ \eqref{D} obeys \eqref{dl9} and, so
\be
\left |\intr |u|(u\cdot \na\beta) dx\right| \le \fr{\nu}{6}D + C\nu^{-1} r(t)^{-2} \|u\|_{L^3}^5 .
\la{betaintbo}
\ee
In view of \eqref{intl34} and the assumption $\int_0^T r(t)^{-4}dt <\infty$, it follows  that the factor
$r^{-2}\|u\|_{L^3}^2$ multiplying $\|u\|_{L^3}^3$ is time integrable a priori,
\be
r(t)^{-2}\|u(t)\|_{L^3}^2 \in L^1([0,T]),
\la{rl3b}
\ee
and thus this term leads to an explicit  uniform bound on $\|u\|_{L^3}$, in terms of the initial data and $\int_0^T r^{-4}dt$.
We integrate by parts in the term
\be
\left |\intr |u| (u\cdot \na\pi) dx\right | = \left |\intr \pi (u\cdot \na |u|) dx\right | \le \int_{|u|\le U} |\pi| |u| |\na |u|| dx +
\int_{|u|\ge U} |\pi| |u| |\na |u|| dx
\la{intpisplit}
\ee
where $U$ is a large time independent constant, at our disposal. We estimate the first term using \eqref{lqbounds}:
\be
\int_{|u|\le U} |\pi| |u| |\na |u|| dx \le C_2 U^{\fr{1}{2}}\sqrt{D} \|u\|_{L^4}^2  \le \fr{\nu}{6} D + C U^2\nu^{-3} \|u\|_{L^3}^5.
\la{firstpibo}
\ee
where we used interpolation (used also in \eqref{pl2interp}) and \eqref{dl9}. We argue like in the proof of Theorem \eqref{pl32thm}, invoking
\eqref{intl34}, which implies that the term $U^2 \|u\|_{L^3}^2$ multiplying $\|u\|_{L^3}^3$ in the right hand side of \eqref{firstpibo} is time integrable with an explicit a priori bound, and as such it leads via Grownwalll to an explicit bound on $\|u\|_{L^3}$.

Finally, taking $U$ large enough so that the set $B_U(t)$ of \eqref{bu} has small measure as in \eqref{bumeas}, using \eqref{pis2b}, proceeding in the same manner as for \eqref{twoall}, and using the assumption \eqref{S2cond}, we have
\be
\int_{|u|\ge U} |\pi| |u| |\na |u|| dx \le CD\left(\int_{B_U}S_2(\cdot, r(t))^{\fr{3}{2}}dx\right)^{\fr{1}{3}} \le \fr{\nu}{6}D.
\la{secondpib}
\ee
This term is absorbed in the remaining dissipative term, ending the proof.
\end{proof}
As in the case of the condition regarding the finite uniform integrability of the pressure \eqref{unifint},  the structure function finite integrability condition \eqref{S2cond} leads to easily accessible bounds.
\beg{thm}\la{S2quanta} Let $q\ge 4$. There exists a constant $C$ depending on $q$ such that if \eqref{S2cond} holds on $[0,T]$ then 

\be
u\in L^{\infty}(0,T; L^q(\Rr^3))
\la{ulqS2b}
\ee
holds. More precisely, we have the single exponential bound
\be
\|u(\cdot,t)\|_{L^q} \le \|u_0\|_{L^q} \exp{\left (\fr{C t\|u_0\|^2_{L^2} }{\nu\delta} + \nu^{-\fr{3}{2}}\|u_0\|_{L^2}^2\Gamma (t)\right)}
\la{lqbound}
\ee
where $\delta$ is the constant in \eqref{S2cond} and
\be 
\Gamma(t) = \sqrt{\int_0^t r^{-4}(s)ds}
\la{Gamma}
\ee
is bounded on $[0,T]$ by assumption.
\end{thm}
\beg{proof}
The proof follows closely the proof of Theorem \ref{presl32quanta}.  We use the evolution of the $L^q$ norm \eqref{lrev}
where we changed $r$ to $q$ because $r$ has now a different meaning. We split at each time $p = \beta(\cdot, r) + \pi(\cdot,r)$ using $r=r(t)$.
We bound the term
\be
\intr \beta u\cdot\na |u|^{q-2}dx \le \|\beta\|_{L^q}\|u\|_{L^q}^{\fr{q-2}{2}}\sqrt D\le C\|\beta\|_{L^{\infty}}^{\fr{1}{2}}\|u\|_{L^q}^{\fr{q}{2}}\sqrt D
\la{betabounds}
\ee
where we used H\"{older} with exponents $q$, $2$ and $\fr{2q}{q-2}$, interpolated $\|\beta\|_{L^{q}}\le \|\beta\|_{L^{\infty}}^{\fr{1}{2}} \|\beta\|_{L^{\fr{q}{2}}}^{\fr{1}{2}}$ and    used the bound $\|\beta\|_{L^{\fr{q}{2}}} \le C\|u\|_{L^q}^2$.
Hiding $\sqrt{D}$ in $\fr{1}{3}\nu D$, this term leads to a growth factor 
\be
\nu^{-1}\int_0^t \|\beta(\cdot, r(s))\|_{L^{\infty}}ds \le C\nu^{-\fr{3}{2}}\|u_0\|_{L^2}^2 \Gamma(t)
\la{betiftyb}
\ee
where we used \eqref{betifty} and the Navier-Stokes energy inequality.  We treat the terms involving $\pi$ in exactly the same manner as we treated $p$ in the proof of Theorem \ref{presl32quanta}.  We omit further details.
\end{proof}
\beg{rem} Theorem \ref{S2quanta} is stronger than Theorem \ref{S2reg} for strong solutions, which have initial data in $H^1$. 
\end{rem}
\beg{rem} A particularly significant small scale $r =\eta $  is given by classical turbulence theory, where the Kolmogorov dissipation wave number $k_d$, inverse of the viscous dissipation scale $\eta$ is  given by
\be
k_d = \eta^{-1} = {\nu}^{-\fr{3}{4}} \epsilon^{\fr{1}{4}}  = \nu^{-\fr{1}{2}} (\langle |\na u (t)|^2\rangle)^{\fr{1}{4}}.
\la{kd}
\ee
We note that it is a priori time integrable to power $4$.  
\end{rem}
\subsection{A nearly selfsimilar example.} 
Theorem \ref{S2reg}  (or rather, its proof) applies to functions which have small translation increments in $L^3$. We consider 
\be
u(x,t) = V + \text{smooth}
\la{uselfsim}
\ee
where the leading term $V$ satisfies
\be
\|V(y+\cdot)- V(\cdot)\|_{L^3} \le W(t) |y|^s
\la{Vprofile}
\ee
for some $s>0$. 
In order to have nondimensional quantities, we write 
\be
W = UL^{1-s}.
\la{WUL}
\ee
The typical example is of the form
\be
V(x,t) = U(t) P\left( \fr{x}{L(t)}\right)
\la{vup}
\ee
where $P$ is time independent and $\|P(\cdot +z)-P(\cdot)\|_{L^3}\le C|z|^s$. We note that this condition is satisfied by many functions with slow decay which are not in $L^3(\Rr^3)$ or even in $L^2(\Rr^3)$, such 
as $P(z) = (1+|z|)^{-\beta}$, $\beta>0$. Of course, the condition is also satisfied on $B^s_{3,\infty}(\Rr^3)$. 

We have that \eqref{Vprofile} reads
\be
\|V(y+\cdot)- V(\cdot)\|_{L^3} \le UL \left( \fr{|y|}{L}\right)^s,
\la{VprofileUL}
\ee
and define 
\be
\fr{UL}{\nu} = Re(V).
\la{scaling}
\ee
We take $\epsilon>0$ such that $s>\fr{1}{2}\epsilon$, and, writing $|y|^{-3} = |y|^{-1+\epsilon}|y|^{-2-\epsilon}$ we use a H\"{o}lder inequality with exponents $3, \fr{3}{2}$  to bound
\be
S_2(x,r)^{\fr{3}{2}}\le C{\epsilon}^{-\fr{1}{2}} {r}^{\fr{3\epsilon}{2}}\int_{|y|\le {r}}|y|^{-3-\fr{3\epsilon}{2}}\left |V(y+x)- V(x)\right|^3dy\\
\la{s_2}
\ee
Integrating $dx$ on $A$ and switching the order of integration we deduce
\be
\ba
\int_A S_2(x,r)^{\fr{3}{2}} dx \le C{\epsilon}^{-\fr{1}{2}} {r}^{\fr{3\epsilon}{2}} \int_{|y|\le {r}}|y|^{-3-\fr{3\epsilon}{2}}\int_{A}\left |V(y+x)- V(x)\right|^3dxdy\\
\le C{\epsilon}^{-\fr{1}{2}}  W^3 r^{3s} =   C{\epsilon}^{-\fr{1}{2}}  (UL)^3\left(\fr{r}{L}\right)^{3s}.
\ea
\la{s2selfb}
\ee
Assuming a bound on the Reynolds number of the profile,
\be
Re(V) \le R
\la{revr}
\ee
and fixing $\epsilon <2s$, we have
\be
\int_A S_2(x,r)^{\fr{3}{2}} dx \le C_s\left(\fr{r}{L}\right)^{3s} (Re(V)^3\nu^3\le C_s \left(\fr{r}{L}\right)^{3s} R^3\nu^3 \la{s2selfbo}
 \ee
The condition \eqref{S2cond} is satisfied if 
\be
\left(\fr{r}{L}\right)^{s}R \le C_s^{-\fr{1}{3}} (2C)^{-1}.
\la{condr}
\ee
\beg{rem}
The condition \eqref{condr} shows that for self-similar profiles with time dependent collapsing inner scale $L$, the condition is satisfied choosing $r$ small compared to the collapsing scale $L$.  Regularity follows if $L^{-4}(t)$ is time integrable.
In particular, if the leading term $V$ is given by \eqref{vup} and $U(t)L(t) \le R\nu$, then regularity follows. The proof of
this fact follows verbatim the proof of Theorem \eqref{S2reg}  including the estimate \eqref{firstpibo}. In that estimate now $U$ is time dependent and it is bounded above by $RL(t)^{-1}$. The  term $U^{-2}\|u\|_{L^3}^2$ is still time integrable and that is why the result continues to hold.
\end{rem}

\subsection{A Dini Condition}
\beg{thm}\la{dini}
Assume that $u$ satisfies
\be
\|\delta_y u\|_{L^3} \le m(|y|)
\la{besm}
\ee
where $\delta_y u (x,t) = u(x+y, t)-u(x,t)$, and where $0\le m$ is a time independent function satisfying
\be
\int_0^1 m^2(\rho)\fr{d\rho}{\rho}  <\infty.
\la{mcond}
\ee
Then $u$ satisfies \eqref{S2cond} with $r$ time independent, and consequently, smooth solutions of Navier-Stokes equations obeying \eqref{besm} with \eqref{mcond} on $[0,T)$ obey $u\in L^{\infty}(0,T; L^3(\Rr^3))$ with explicit bounds depending only on $m, \nu, T, \|u_0\|_{L^3}, \|u_0\|_{L^2}$.
\end{thm}
\beg{proof}
The proof follows from the fact that
\be
\|S_2(\cdot, r)\|_{L^{\fr{3}{2}}} \le C \int_{|y|\le 2r}\|\delta_y u\|^2_{L^3}\fr{dy}{|y|^3}.
\la{S2bes}
\ee
This inequality is proved by duality, integrating $S_2$ against a test function in $L^3$
\be
\intr S_2 \phi dx = \fr{1}{4\pi}\int_{|y|\le 2r}\fr{dy}{|y|^3}\intr \phi(x) |\delta_y u(x)|^2 dx\le \fr{1}{4\pi}\|\phi\|_{L^3}\int_{|y|\le 2r}\|\delta_y u\|_{L^3}^2 \fr{dy}{|y|^3}.
\ee
From \eqref{S2bes} and the assumed Dini condition, we deduce that the $\|S_2(\cdot, r)\|_{L^{\fr{3}{2}}} \le\left(\fr{\nu}{C}\right)^2$
if $r$ is chosen small enough so that
\be
\int_0^{2r} m^2(\rho)\fr{d\rho}{\rho} \le \left(\fr{\nu}{C}\right)^2.
\la{mr}
\ee
\end{proof}
\beg{rem} Clearly $m(r) \sim \log^{-\alpha}(r^{-1})$ with $\alpha>\fr{1}{2}$ is sufficient. As we remarked before, the smallness of the $L^3$ increment does not imply that the function needs to be in $L^3$. We also remark that $m$ can be allowed to depend on time, if $m(r)r^{-1}$ is uniformly integrable on $[0,1]$, or more generally, if denoting 
\be
I_{m(t)}(r) = \int_0^{2r} m^2(\rho,t)\fr{d\rho}{\rho}
\la{imrt}
\ee
we have that the preimage of $\left(\fr{\nu}{C}\right)^2$ under $I_{m(t)}$,  that is $r(t) = I_{m(t)}^{-1}\left(\left(\fr{\nu}{C}\right)^2\right)$, obeys $\int_0^T r(t)^{-4}dt<\infty$.
\end{rem}
\subsection{Multifractal intermittent scenario}\la{multi}
We consider the region 
$B_U(t) = \{x\left |\right.\; |u(x,t)| \ge U\}$ defined before  in \eqref{bu}. We take
$U$ time independent. We introduce a time independent length scale $L>0$, and we require
\be
U^2  \ge L^{-3}\intr |u|^2dx
\la{UL}
\ee
so that
\be
|B_U| \le L^3.
\la{bul}
\ee
We assume that  the velocity increments 
\be
s_2(x,r)  = \fint_{|y| =r} |u(x+y)-u(x)|^2dS(y)
\la{s2}
\ee
obey bounds
\be
s_2(x,r) \le G^2 \left(\fr{r}{L}\right)^{2\alpha(x)}
\la{s2alpha}
\ee
with $G>0$ constant (with units of velocity), $L>0$ as above, constant, (with units of length)  and with $0<\alpha(x)\le 1$. This upper bound is assumed to hold a.e. in  $x\in B_U(t)$ and for all $0<r<r_0$, where $0< r_0<L$ is a fixed positive constant. Because
\be
S_2(x,r) = \int_0^{2r} s_2(x,\rho) \fr{d\rho}{\rho},
\la{S2s2}
\ee
we have that
\be
S_2(x, r_0) \le C G^2\fr{1}{\alpha(x)} \left (\fr{r_0}{L}\right)^{2\alpha(x)}
\la{S2s2b}
\ee
holds  a.e in  $x\in B_U$. 
In multifractal turbulent intermittent scenarios, it is assumed that there is a spectrum of near-singularities of H\"{o}lder exponent $h$ and that these are achieved on sets $\Sigma_h$ of dimension $d(h)\le 3$ which occur randomly with  probability $d\mu(h)$.

The dimension $d(h)$ is implemented in the following manner. We take a region $V_h$ around $\Sigma_h$ 
and partition it in small disjoint cubes of size $\rho$ with  $\rho\le r_0$.  This region is a "collar"of cross-section size $\rho$ around the set $\Sigma_h \cap B_U$. The multifractal assumption is that the number of such cubes of $V_h$  is of the order
$N_h(\rho) = \left(\fr{\rho}{L}\right)^{-d(h)}$. Assuming $\alpha(x)\ge h$ to hold on each such cube, we have from
\eqref{S2s2b}, on each cube
\be
S_2(x,r_0) \le CG^2 h^{-1} \left(\fr{r_0}{L}\right)^{2h}.
\la{S2h}
\ee
Writing the volume of the cube as $L^3(\fr{\rho}{L})^3$, we have
\be
\int_{B_U\cap V_h} S_2(x, r_0)^{\fr{3}{2}}dx \le C(G L)^3 \fr{1}{h^{\fr{3}{2}}}\left(\fr{r_0}{L}\right)^{3h} \left(\fr{\rho}{L}\right)^{3-d(h)} \le C(GL)^3 h^{-\fr{3}{2}} \left(\fr{r_0}{L}\right)^{3-d(h) + 3h}.
\la{intsh}
\ee
Above we used $\rho \le {r_0} $.
Summing in $h$, remembering the frequency, we obtain
\be
\int_{B_U\cap (\cup_h V_h)}S_2^{\fr{3}{2}}(x,r_0)dx \le C(GL)^3 \int_0^1  h^{-\fr{3}{2}}\left (\fr{r_0}{L}\right)^{3-d(h) + 3h}d\mu(h)
\la{intas1}
\ee
In the multifractal formalism, the structure function exponents are defined by
\be
\zeta_p = \inf_h (3-d(h) +ph).
\la{zetap}
\ee
The inequality \eqref{intas1}  above implies
\be
\int_{B_U\cap (\cup_h V_h)}S_2(x, r_0)^{\fr{3}{2}} dx \le C_\mu (GL)^3 \left(\fr{r_0}{L}\right)^{\zeta_3}
\la{intas2}
\ee 
where
\be
C_\mu = C\int_0^1 h^{-\fr{3}{2}}d\mu(h)
\la{cmu}
\ee
is assumed to be finite.  Introducing the Reynolds number based on $G$,
\be
R_G = \fr{GL}{\nu}
\la{RG}
\ee
and recalling that $B_U \setminus \cup_h V_h$ was assumed to have measure zero, we have
\be
\int_{B_U} S_2(x,r_0)^{\fr{3}{2}} dx \le \left[C_\mu R_G^3 \left(\fr{r_0}{L}\right)^{\zeta_3}\right]
\nu^3
\la{intas3}
\ee
we see that the condition \eqref{S2cond} is satisfied for $B_U$  if
\be
R_G^3 \left(\fr{r_0}{L}\right)^{\zeta_3} \le (C^3C_{\mu})^{-1}.
\la{RGcond}
\ee
In classical turbulence theory $\zeta_3=1$. If $\zeta_3>0$, under the above scenario, it is enough to have $\fr{r_0}{L}$ small enough in order to deduce that no singularities in finite time can occur. 

\subsection{Time dependent regions of interest.}\la{time}

As me noted before, the finite uniform integrability of condition \eqref{S2cond} is not needed, all we need is control of $S_2$ on certain small sets of interest. We consider the set  $B_{U,G}(t)  = \{x\left | \right.\;  |u(x,t)| \ge U, \; \text{and}\; |\na u(x,t)| \ge G\}$ defined in \eqref{bug}. We note that
\be
\left |B_{U,G}(t)\right | \le C\min\{ U^{-2}\|u_0\|_{L^2}^2;  G^{-1}\|\omega_0\|_{L^1}\}.
\la{measB}
\ee
where $\omega = \na\times u$. The first term in the inequality follows from the Markov-Chebyshev inequality and the fact that the $L^2$ norms of solutions of Navier-Stokes equations are non-increasing in time. The second term follows from the fact that the map $\omega\mapsto \na u$ is weak type 1, that is from $G|\{x \left|\right. \; |\na u|\ge G\}\le C\|\omega\|_{L^1}$, and the fact that the $L^1$ norm of vorticity of solutions of Navier-Stokes equations is non-increasing in time \cite{carea}.

\beg{thm}\la{S2onB} Let $U(t)$, $G(t)$ and $r(t)$ be positive numbers such that
\be
\int_0^T (r(t)^{-4} + U(t)^4 + G(t))dt <\infty.
\la{rUGb}
\ee
Consider the set 
\be
B(t) =\{x\left |\right. |u(x,t)|\ge U \; \text{and} \; |\na u(x,t)| \ge G\}.
\la{Bt}
\ee
There exists an absolute constant $C$ such that, if
\be
\int_{|y|\le r(t)}  \left (\int_{B(t)} |\delta_y u(x,t)|^3dx\right)^{\fr{2}{3}}\fr{dy}{|y|^3}  \le \left( \fr{\nu}{C}\right)^2
\la{BS2cond}
\ee
then the smooth solution of Navier-Stokes equations obeys $u\in L^{\infty}(0,T; L^3(\Rr^3))$ with explicit bounds depending only on $\nu, T, \|u_0\|_{L^3}, \|u_0\|_{L^2}$.
\end{thm}
\beg{proof}
We follow the proof of Theorem \ref{S2reg}. The $\beta$ term is estimated as in \eqref{betaintbo}. The contribution of 
the term involving $\pi$ from the region $|u|\le U$ is estimated as in \eqref{firstpibo}, noting that the term $U^2\|u\|_{L^3}^2$ is time integrable in view of \eqref{rUGb}. A new term is
\be
\int_{|u| \ge U, |\na u|\le G} |\pi| |u| |\na |u|| dx \le CG\|u\|_{L^3}^3,
\la{piUGb}
\ee
and, in view of \eqref{rUGb} it leads via Grownwalll to an explicit bound on $\|u\|_{L^3}$. We are left with
\be
\int_{B(t)} |\pi| |u| |\na |u|| dx \le C D\left(\int_{B(t)} S_2^{\fr{3}{2}}(x, r(t))dx\right)^{\fr{1}{3}}
\la{intB}
\ee
Now we use 
\be
\left(\int_B S_2^{\fr{3}{2}} dx \right)^{\fr{2}{3}} \le \fr{1}{4\pi}\int_{|y| \le r} \fr{dy}{|y|^3}\int_B |\delta_y u(x)|^2dx
\la{normy}
\ee
proved by duality, testing against arbitrary $L^3(B)$ functions. The assumption \eqref{BS2cond} implies
\be
\int_{B(t)} |\pi| |u| |\na |u|| dx \le \fr{\nu D}{6}
\la{intBb}
\ee
and concludes the proof.
\end{proof}

\beg{thm}\la{BS2quanta} Let $q\ge 4$, and assume \eqref{BS2cond} where the functions $U(t)$, $r(t)$ and $G(t)$
obey
\be
\int_0^T (U^2(t) + r^{-4}(t) + G(t))dt <\infty.
\la{UrGb}
\ee
Then we have the single exponential bound
\be
\|u(\cdot, t)\|_{L^q}\le \|u_0\|_{L^q}\exp{\left( C\nu^{-1}\int_0^t U^2ds + C\int_0^t G(s)ds + C\nu^{-\fr{3}{2}}\|u_0\|_{L^2}^2\sqrt{\int_0^t r^{-4}(s)ds}\right )}
\la{lqBS2b}
\ee
\end{thm}
\beg{proof} 
We start as in the proof of Theorem \ref{S2quanta} by splitting $p = \beta + \pi$ in the estimate the evolution of the $L^q$ norm of $u$, and deduce the bound \eqref{betabounds} leading to the exponential growth factor \eqref{betiftyb}.
We are left o estimate the contribution of $\pi$, that is
\be
I = \intr \pi u\cdot \na |u|^{q-2}dx.
\la{picontro}
\ee
We bound the integral 
\be
\ba
|I| \le  I_1 + I_2 + I_B \\
= (q-2)\left(\int_{|u|\le U} |\pi ||u|^{q-2}|\na u| dx + \int_{|\na u|\le G}  |\pi ||u|^{q-2}|\na u| dx + \int_{B(t)} |\pi ||u|^{q-2}|\na u| dx\right).
\ea
\la{I123}
\ee
We bound $I_1$ like in \eqref{rle4},
\be
I_1 \le CU\|u\|_{L^q}^{\fr{q}{2}}\sqrt{D}
\la{I1b}
\ee
where we use the fact that $\|\pi (\cdot, r)\|_{L^{\fr{q}{2}}}\le C \|u\|_{L^q}^2$ holds with $C$ an absolute constant, independent of $r$, and 
\be
\left(\int_{|u|\le U}|u|^{\fr{q(q-2)}{q-4}}dx\right)^{\fr{q-4}{2q}} \le U\|u\|_{L^q}^{\fr{q}{2}-2}.
\la{uUqr}
\ee
 The bound \eqref{I1b} is valid for $q=4$ as well, we just take $|u|\le U$ outside the integral and use $L^2-L^2$ bounds.
The term $I_2$ is bound directly
\be
|I_2| \le CG\|u\|_{L^q}^q
\la{I2b}
\ee 
The last term is smaller than the dissipation, using the arguments similar to the ones leading to \eqref{twotwor}.
 We omit further details.
\end{proof}

\beg{rem} The condition $U\in L^2(0,T)$ appearing in \eqref{UrGb} is better than the condition $U\in L^4(0,T)$ of 
\eqref{rUGb} of Theorem \ref{S2onB}. That is just because in that theorem the desire was to bound the $L^3$ norm in terms solely of itself. Theorem \ref{BS2quanta} is strictly stronger that Theorem \ref{S2onB} (it implies it for strong solutions), by bounding first the $L^4$ norm of the solution, and then returning to the proof of the bound of the $L^3$ norm.
\end{rem}

{\bf{Acknowledgment.}} Research partially supported by NSF grant DMS-2106528.

\end{document}